\documentclass{article}
\usepackage{arxiv}

\usepackage[utf8]{inputenc} 
\usepackage[T1]{fontenc}    
\usepackage{amsmath}
\usepackage{authblk}
\usepackage{amsfonts,amssymb}
\usepackage{amsthm}
\usepackage{xcolor}

 \usepackage[pagewise]{lineno}

\newcommand{\R}{\mathbb{R}}

\newtheorem{theorem}{Theorem}[section]
\newtheorem{proposition}{Proposition}[section]

\usepackage{scalerel,stackengine}
\stackMath
\newcommand\reallywidehat[1]{%
\savestack{\tmpbox}{\stretchto{%
  \scaleto{%
    \scalerel*[\widthof{\ensuremath{#1}}]{\kern-.6pt\bigwedge\kern-.6pt}%
    {\rule[-\textheight/2]{1ex}{\textheight}}
  }{\textheight}%
}{0.5ex}}%
\stackon[1pt]{#1}{\tmpbox}%
}
\usepackage{scalerel}[2014/03/10]
\usepackage{stackengine}

\newcommand{\bb}{\begin{equation}}
\newcommand{\ee}{\end{equation}}
\newcommand{\ba}{\begin{array}}
\newcommand{\ea}{\end{array}}
\newcommand{\f}{\frac}
\usepackage[all]{xy}

\numberwithin{equation}{section}

\newcommand{\sobolev}[2]{\Vert #1\Vert_{H^{#2}}}

\title{Analytic regularity of global solutions for the $b$-equation}

    \author[1,2]{\Large \bf Priscila Leal da Silva}
    \affil[1]{\large Department of Mathematical Sciences, Loughborough University, United Kingdom}
    \affil[2]{\large Centro de Matem\'atica, Computa\c{c}\~ao e Cogni\c{c}\~ao, Universidade Federal do ABC, Brazil}
    \arxiv{AAAA}

\begin{document}

\email{P.Leal-Da-Silva@lboro.ac.uk}

\maketitle

\begin{abstract}
In this paper, we delve into the $b$-family of equations and explore regularity properties of its global solutions. Our findings reveal that, irrespective of the real choice of the constitutive parameter, when the initial datum is confined to an analytic Gevrey function the resulting global solution is analytic in both temporal and spatial variables.
\end{abstract}



\section{Introduction}
This paper is concerned with the $b$-family of equations
\begin{align}\label{eq1.1}
    \begin{aligned}
        u_t= - uu_x - \partial_x\Lambda^{-2}\left(\f{b}{2}u^2 + \f{3-b}{2}u_x^2\right),
    \end{aligned}
\end{align}
where $b$ is a real number and $\Lambda^{-2}$ denotes the inverse of the Helmholtz operator $\Lambda^2 =1-\partial_x^2$ acting on Sobolev spaces, and the analytic regularity of its global solutions. Equation \eqref{eq1.1} was first considered in \cite{DHH,HS1,HS2} with applications in hydrodynamics for $b\neq -1$ and its most prominent members are the Camassa-Holm \cite{CH} equation ($b=2$) and the Degasperis-Procesi \cite{DP} equation ($b=3$), the only two cases where there exists an infinite number of linearly independent conserved vectors.

The classical theory for global well-posedness of non-local equations of the type \eqref{eq1.1} can be attributed to Constantin and Escher \cite{CE}. Their work demonstrated the global well-posedness of the Camassa-Holm equation in $C([0,\infty), H^3(\R))$ for an initial datum $u_0(x):=u(0,x) \in H^3(\R)$, subject to the condition that the initial momentum $m_0(x):=u_0(x)-u_0''(x)$ does not change sign. Their proof heavily relied on the conservation of energy expressed in terms of the usual norm in $H^1(\R)$.

For the Degasperis-Procesi equation, Liu and Yin \cite{LY} provided a different proof that does not rely on conservation of energy, as its solutions do not conserve the $H^1(\R)$-norm. Instead, they accomplished a similar result for $u_0\in H^s(\R)$, where $s>3/2$, by estimating the $L^{\infty}(\R)$ norm. 

For arbitrary choices of $b$, however, the conservation of energy does not hold. In fact, one of the few conserved quantities is the $L^1(\R)$-norm of $m_0(x)$ when it is either non-negative or non-positive. However, relying on this norm is often insufficient for establishing well-posedness and extending local solutions of \eqref{eq1.1} to a global scope, presenting a significant challenge.

In the papers \cite{GLT,LY}, the unique global strong solutions of equation \eqref{eq1.1} were proven to exist for an initial datum $u_0(x)\in H^s(\R)$, with $s>3/2$, but only for the case when $b>0$. Global well-posedness for the scenario where $b=0$ was obtained in \cite{dSF}. Moreover, in \cite{EY}, Escher and Yin extended the global solutions to the case where $b=-1/2n$, where $n$ is any positive integer. Recent contributions by Freire in \cite{F} further extended the analysis of global solutions to encompass any $b\in\R$, thus providing a complete and comprehensive understanding of the existence and uniqueness of global solutions of \eqref{eq1.1}.

Building upon \cite{Eu} and incorporating the recent findings from \cite{F}, we investigate \eqref{eq1.1} with initial data belonging to the Gevrey class $G^{\sigma,s}(\R)$ of functions $f\in L^2(\R)$ such that the norm
\begin{align}
\Vert f\Vert_{G^{\sigma,s}}=\left(\int_{\R}e^{2\sigma |\xi|}(1+\xi^2)^s|\hat{f}(\xi)|^2d\xi\right)^{1/2}
\end{align}
remains finite, with $\sigma>0$ and $s\in\R$. Here, $\hat{f}(\xi)$ denotes the Fourier transform of $f$, and for $s\geq 0$, the continuous embedding $G^{\sigma,s}(\R)\hookrightarrow H^{\infty}(\R)$ holds, where $H^{\infty}(\R) := \bigcap\limits_{s\geq 0}H^{s}(\R)$.

In \cite{Eu}, the author showed that when the initial data lies in the Gevrey space $G^{1,s}(\R)$ of analytic functions, with $s>5/2$, the solution of \eqref{eq1.1} with $b=0$ is globally analytic in both variables. The primary objective of this paper is to extend this crucial result to encompass the entire family \eqref{eq1.1} for $b\in \mathbb{R}$ by establishing the following theorem.

\begin{theorem}\label{thr1}
    Given $u_0(x)\in G^{1,s}(\R)$, for $s>3/2$, such that $m_0(x)\geq 0$ or $m_0(x)\leq 0$, there exists a unique global analytic solution $u\in C^{\omega}([0,\infty)\times \R)$ of \eqref{eq1.1}.
\end{theorem}

For this purpose, we will prove some auxiliary local well-posedness results for analytic initial datum and combine these results with the global well-posedness proved in \cite{F} and the Kato-Masuda machinery \cite{KM} in order to complete the proof.

Theorem \ref{thr1} represents a notable extension of the existing literature. While \cite{BHP1} has established the analytic regularity of solutions for the Camassa-Holm equation, the current theorem transcends the limitations imposed by the (in general) lack of conservation laws and rendering it independent of the specific choice of $b$ in \eqref{eq2.1}. This, in turn, presents the challenge of enhancing the radius of spatial analyticity (as elaborated in Section \ref{sec3}) with a  polynomial lower bound for it.


The paper is structured as follows:
\begin{itemize}
    \item in Section \ref{sec2}, we present the auxiliary spaces and results required to establish (analytic in time) local well-posedness in Gevrey and Himonas-Misiolek \cite{HM} spaces;
    \item Section \ref{sec3} utilizes the Kato-Masuda \cite{KM} machinery to obtain spatial analytic solutions, along with determining the radius of spatial analyticity;
    \item Lastly, in Section \ref{sec4}, we employ space embeddings and the results from Sections \ref{sec2}-\ref{sec3} to finalize the proof of Theorem \ref{thr1}.
\end{itemize}

\section{Function spaces and local well-posedness}\label{sec2}

In this section we present the remaining function spaces needed and also the auxiliary local results. Besides Gevrey spaces, we will consider the Himonas-Misiolek space \cite{HM}
$$E^{\sigma,m}(\R)=\{f\in C^{\infty}(\R); \Vert f\Vert_{E^{\sigma,m}} = \sup\limits_{j\in\mathbb{Z}_+}\frac{\sigma^j(j+1)^2}{j!}\Vert \partial_x^jf\Vert_{H^{2m}}<\infty\}$$
for $\sigma>0$ and $m$ being a positive integer.

Let $F(u) = - uu_x - \partial_x\Lambda^{-2}\left(\frac{b}{2}u^2 + \frac{3-b}{2}u_x^2\right)$ represent the right-hand side of equation \eqref{eq1.1} and consider initial data $u_0\in G^{\sigma_0,s}(\R)$ and $v_0\in E^{\sigma_0,m}(\R),$ where $s>3/2$, $m\geq 2$ and $0<\sigma_0\leq 1$. If $\sigma \in (0,\sigma_0)$, $u_i\in G^{\sigma,s}(\R)$ and $v_i\in E^{\sigma,m}(\R)$, for $i=1,2$, are such that
\begin{align*}
    \Vert u_i-u_0\Vert_{G^{\sigma,m}}<\Vert u_0\Vert_{G^{\sigma_0,s}},\quad  \Vert v_i-v_0\Vert_{E^{\sigma,m}}<\Vert v_0\Vert_{E^{\sigma_0,m}}, \forall i,
\end{align*}
then by employing a standard argument with the use of algebra properties of Gevrey (Lemma 3 in \cite{BHP}) and Himonas-Misiolek (Theorem 2.1 in \cite{HM}) spaces, along with their corresponding derivatives estimates, see Lemma 2.4 in \cite{HM}, Lemma 2 in \cite{BHP} and Lemma 2.2 in \cite{Eu}, we can derive the following estimates:
\begin{align*}\label{eq2.1}
    \begin{aligned}
        &\Vert F(u_0)\Vert_{G^{\sigma_0,s}}\leq \f{M_1(\Vert u_0\Vert_{G^{\sigma_0,s}},b)}{\sigma-\sigma_0},\quad \Vert F(v_0)\Vert_{E^{\sigma_0,s}}\leq \f{M_2(\Vert v_0\Vert_{E^{\sigma_0,m}},b)}{\sigma-\sigma_0}\Vert v_0\Vert_{E^{\sigma,m}}\\
        &\Vert F(u_1)-F(u_2)\Vert_{G^{\sigma_0,s}}\leq \f{L_1(\Vert u_0\Vert_{G^{\sigma_0,s}},b)}{\sigma-\sigma_0}\Vert u_1-u_2\Vert_{G^{\sigma,s}},\\
        &\Vert F(v_1)-F(v_2)\Vert_{E^{\sigma_0,s}}\leq \f{L_2(\Vert v_0\Vert_{E^{\sigma_0,m}},b)}{\sigma-\sigma_0}\Vert v_1-v_2\Vert_{E^{\sigma,m}},
    \end{aligned}
\end{align*}
for some positive $M_1,M_2,L_1$ and $L_2$ depending only on $b$, the respective initial data and parameters of the spaces under consideration. From the Autonomous Ovsyanikov Theorem (see Theorem 1 in \cite{BHP}), we prove the following two results:
\begin{proposition}\label{prop2}
     Given $u_0\in G^{\sigma_0,s}(\R)$, with $s>3/2$ and $\sigma_0\in(0,1]$, there exists $T=T(\Vert u_0\Vert_{G^{\sigma_0,s}})>0$ such that for every $\sigma\in (0,\sigma_0)$ the initial value problem \eqref{eq1.1} has a unique solution $u\in C^{\omega}([0,T(1-\sigma)), G^{\sigma,s}(\R))$.
\end{proposition}

\begin{proposition}\label{prop3}
    Given $u_0\in E^{\sigma_0,m}(\R)$, with $\sigma_0\in (0,1]$ and $m\geq 2$, then there exists $\epsilon=\epsilon(\Vert u_0\Vert_{E^{\sigma_0,m}})>0$ and a unique solution $u\in C^{\omega}([0,\epsilon],E^{\sigma,m}(\R))$  of \eqref{eq1.1} for any $\sigma \in (0,\sigma_0)$.
\end{proposition}

While Proposition \ref{prop3} will be used for general values of $\sigma_0$, our only interest in Proposition \ref{prop2} is regarding the case $\sigma_0=1$ as the corresponding Gevrey spaces consist of analytic functions. Moreover, observe that in Proposition \ref{prop3} we can extend the existence interval to $[-\epsilon,\epsilon]$ by making a reflection of $t$ and $x$.

\section{Radius of spatial analyticity}\label{sec3}

In this section we will prove spatial analyticity of the solution by making use of the Kato-Masuda machinery, see Theorem 1 in \cite{HM}. But before we need to build the settings in which we will work.

In \cite{F}, the author showed that for any choice of $b\in\R$, given an initial datum $u_0(x)\in H^{n+2}(\R)$, for any positive integer $n\in \mathbb{N}_0:=\{0,1,2,3,\dots\}$, such that $m_0(x)\geq 0$ or $m_0(x)\leq 0$, then there exists a unique global solution $u\in C([0,\infty),H^{n+2}(\R))$ of \eqref{eq1.1}. On the one hand, since $G^{1,s}(\R)\hookrightarrow H^{\infty}(\R)$, then any initial datum $u_0\in G^{1,s}(\R)$, with $s>3/2$ is in $H^{n+2}(\R)$ for any $n$, guaranteeing the existence of global solutions $u\in C([0,\infty),\bigcap\limits_{n\in\mathbb{N}_0}H^{n+2}(\R))$. On the other hand, fixed $n\in \mathbb{Z}$, we have $H^{n+1}(\R)\hookrightarrow H^s(\R) \hookrightarrow H^{n}(\R)$ for any real number $s\in[n,n+1]$. This means that if $u(t,\cdot)\in \bigcap\limits_{n\in\mathbb{N}_0}H^{n}(\R)$, then $u(t,\cdot)\in H^{\infty}(\R)$.

In the particular case $n=1,$ the solution $u(t,\cdot)$ belongs to $H^3(\R)\subset H^2(\R)\subset H^1(\R)\subset L^2(\R)$ and then we conclude that $u(t,\cdot)\in \bigcap\limits_{n\in\mathbb{N}_0}H^{n}(\R)$, leading to $u(t,\cdot)\in H^{\infty}(\R)$. This summarises the following result.

\begin{proposition}\label{prop4}
    Given an initial datum $u_0(x)\in G^{1,s}(\R)$, with $s>3/2$, such that $m_0(x)\in L^1(\R)\cap H^1(\R)$ is either non-negative or non-positive, then there exist a unique solution $u\in C([0,\infty),H^{\infty}(\R))$ of \eqref{eq1.1}.
\end{proposition}

We now proceed with the Kato-Masuda theorem. Given $u_0\in G^{1,s}(\R)$, with $s>3/2$, let $u$ be the solution whose existence is guaranteed by Proposition \ref{prop4} and define $X=H^{m+2}(\R), Z=H^{m+5}(\R)$, with $m\geq 2$, so that $F:Z\rightarrow X$ is continuous. For $T>0$ fixed, let $\mu = 1+\max\{\Vert u\Vert_{H^2}, t\in[0,T]\}$ and $O=\{v\in Z; \Vert v\Vert_{H^2}<\mu\}$.

For $m\geq 0$ and $\sigma\in\R$ to be better specified, let
\begin{align}
    \Phi_{\sigma,m}(v) = \frac{1}{2}\Vert v\Vert_{\sigma,2,m}^2:=\frac{1}{2}\sum\limits_{j=0}^m\frac{1}{(j!)^2}e^{2\sigma j} \sobolev{\partial_x^ju}{2}^2,
\end{align}
for $v\in O$. With this norm, we define the Kato-Masuda space $A(r)$, for $r>0$, as the set of functions that can be analytically extended to a function on a strip of width $r$, endowed with the norm $\Vert v\Vert_{\sigma,2} = \lim\limits_{m\to\infty}\Vert v\Vert_{\sigma,2,m}$ for every $\sigma\in\R$ such that $e^{\sigma}<r$. We have a similar and useful embedding $G^{\sigma,s}(\R)\hookrightarrow A(\sigma)\hookrightarrow H^{\infty}(\R)$ for the Kato-Masuda space which will be considered in the next section.

For $D\Phi_{\sigma}(v)F(v):=\langle F(v)\,,\, D \Phi_{\sigma}(v)\rangle_{H^2}$, where $D$ denotes the Fréchet derivative and $v\in O$, we can use the triangle inequality to write
    \begin{align*}
        \vert D\Phi_{\sigma,m}(v)F(v)\vert \leq& \left\vert\sum\limits_{j=0}^{m}\frac{e^{2\sigma j}}{(j!)^2}\langle\partial_x^j v\,,\,\partial_x^j(vv_x)\rangle_{H^2} \right\vert+ \frac{|b|}{2}\left\vert\sum\limits_{j=0}^{m}\frac{e^{2\sigma j}}{(j!)^2}\langle\partial_x^j v\,,\,\partial_x^{j+1}\Lambda^{-2}v^2\rangle_{H^2}\right\vert\\
        &+ \frac{|3-b|}{2}\left\vert\sum\limits_{j=0}^{m}\frac{e^{2\sigma j}}{(j!)^2}\langle\partial_x^j v\,,\,\partial_x^{j+1}\Lambda^{-2}v_x^2\rangle_{H^2}\right\vert.
    \end{align*}
By using the estimates (6.14)--(6.16) of \cite{BHP1}, we conclude that for $A(p)= (32+16|b|+64|3-b|)p$ and $B(p,q)=(64+32|b|+256|3-b|)(1+p)q^{1/2}$ we can bound the previous term by 
    \begin{align*}
        \vert D\Phi_{\sigma,m}(v)F(v)\vert \leq& A(\Vert v\Vert_{H^2})\Phi_{\sigma,m}(v) + B(\Vert v\Vert_{H^2},\Phi_{\sigma,m}(v))\partial_{\sigma}\Phi_{\sigma,m}(v).
    \end{align*}
For $K=A(\mu), \beta(p)=Kp$ and $\alpha(p)=B(\mu,p)$, we conclude that 
    \begin{align*}
        \vert D\Phi_{\sigma,m}(v)F(v)\vert \leq& \beta(\Phi_{\sigma,m}(v)) + \alpha(\Phi_{\sigma,m}(v))\partial_{\sigma}\Phi_{\sigma,m}(v),
    \end{align*}
and then Kato-Masuda Theorem (see Theorem 1 in \cite{KM}) yields that for $t\in [0,T]$ the unique solution $u(t)$ belongs to $A(r(t))$, where $r(t)=e^{\sigma(t)}$ and $\sigma(t)=\gamma - \lambda (e^{A(\mu)t/2}-1)$, with $\gamma<0$ fixed and $\lambda$ depending only $b$, the initial datum and $\sigma_0$. This completes the proof of the following result:

\begin{proposition}\label{prop5}
Given $u_0\in G^{1,s}(\R)$, with $s>3/2$, suppose that $m_0\in L^1(\R)\cap H^1(\R)$ does not change sign and let $u\in C([0,\infty);H^{\infty}(\R))$ be the unique solution to the initial value problem of \eqref{eq1.1}. Then for every $T>0$ there exists $r(T)>0$ such that $u\in C([0,T];A(r(T)))$.    
\end{proposition}

We can now proceed with the proof of Theorem \ref{thr1}.

\section{Proof of Theorem \ref{thr1}}\label{sec4}

Once we have proved Propositions \ref{prop2}-\ref{prop5}, the proof of Theorem \ref{thr1} reduces to the space embeddings considered previously \cite{BHP1,Eu}. From now on fix an arbitrary initial datum $u_0(x)\in G^{1,s}(\R)$, with $s>3/2$, and let $u\in C([0,\infty), H^{\infty}(\R))$ be the unique global solution obtained through Propositions \ref{prop4}-\ref{prop5}.

\begin{itemize}
    \item \textbf{Proving that $u\in C^{\omega}([0,T],A(\sigma(T)))$ for some $T>0$ and $\sigma(T)$:}
\end{itemize}    
    Given an initial datum $u_0(x)\in G^{1,s}(\R)$, with $s>3/2$, let $\tilde{u}\in C^{\omega}([0,\tilde{T}(1-\sigma)),G^{\sigma,s}(\R))$ be the unique solution obtained from Proposition \ref{prop2}, where $\sigma\in(0,1)$ and $\tilde{T}>0$.

    By defining $T=\tilde{T}(1-\sigma)/2$, we obtain $\sigma(T) = 1-2T/\tilde{T}$, while the embeddings $G^{\sigma(T),s}(\R)\hookrightarrow A(\sigma(T))\hookrightarrow H^{\infty}(\R)$ yield $\tilde{u}\in C^{\omega}([0,T],A(\sigma(T)))\subset C^{\omega}([0,T],H^{\infty}(\R))$. Consequently, Proposition \ref{prop5} guarantees that $\tilde{u}=u$ for $t\in[0,T]$ and $u\in C^{\omega}([0,T],A(\sigma(T)))$.

    With this, we proved that the solution $u$ is locally analytic in time and is $C^{\infty}(\R)$ in space.

    \begin{itemize}
    \item \textbf{Proving that the analytic lifespan cannot be finite:}
    \end{itemize}  
    In the interval $[0,T]$ constructed in the previous item, let $$T^{\ast} = \sup\{T>0, u\in C^{\omega}([0,T], A(\sigma(T))),\,\,\text{for some}\,\, \sigma(T)>0\}$$
    and assume that $T^{\ast}$ is finite.

    As a consequence of the embedding $A(\sigma(T))\hookrightarrow H^{\infty}(\R)$, Proposition \ref{prop5} establishes that the initial datum $u(T^{\ast})$ lies within $A(r)$. By choosing $\sigma_0<\min\{1, r/e\}$, the inverse of Lemma 5.1 from \cite{BHP1} shows that $u(T^{\ast})$ belongs to $E_{\sigma_0,m}(\R)$ for $m\geq2$ and $\sigma_0\in(0,1]$.
    
    From Proposition \ref{prop4}, we deduce the existence of a unique solution $\tilde{u}\in C^{\omega}([0,\epsilon],E_{\delta,m}(\R))$ for $0<\delta<\sigma_0$, satisfying the initial condition $\tilde{u}(0) = u(T^{\ast})$. Moreover, utilizing the embeddings $E_{\delta,m}(\R)\hookrightarrow A(\delta)\hookrightarrow H^{\infty}(\R)$, we can conclude that $\tilde{u}\in C([0,\epsilon];H^{\infty}(\R))$, implying $\tilde{u}(t) = u(T^{\ast}+t)$ for $t\in[0,\epsilon]$.

    If we let $s:=T^{\ast}+t$, then $u(s) = \tilde{u}(s-T^{\ast}),$ for $s\in[T^{\ast},T^{\ast}+\epsilon],$
    that is,
    $$u \in C^{\omega}([T^{\ast},T^{\ast}+\epsilon];E_{\delta,m}(\R))\subset C^{\omega}([T^{\ast},T^{\ast}+\epsilon];A(\delta)).$$
    Let $T>0$ be a real number satisfying $T^{\ast}-\epsilon< T< T^{\ast}$ so that the solution $u$ belongs to $C^{\omega}([0,T];A(\sigma(T)))$ for some $\sigma(T)>0$.

    By defining $\tilde{\sigma} = \min\{\delta,\sigma(T)\}$, then $$u\in C^{\omega}([0,T];A(\tilde{\sigma})),\quad\text{and}\quad u\in C^{\omega}([T^{\ast}-\epsilon,T^{\ast}+\epsilon];A(\tilde{\sigma})),$$
    which consequently says that $u\in C^{\omega}([0,T^{\ast}+\epsilon];A(\tilde{\sigma}))$ and $T^{\ast}$ cannot be the supremum. As a result of the contradiction, $T^{\ast}$ must be infinite and, for every $T>0$, there exists $r(T)>0$ such that $u\in C^{\omega}([0,T];A(r(T)))$.

\begin{itemize}
    \item \textbf{Conclusion of the proof:}
\end{itemize}
    The proof of Theorem \ref{thr1} is concluded by using a result proved by Barostichi, Himonas and Petronilho in \cite{BHP1}, see page 752, stating that the previous item suffices to prove that $u \in C^{\omega}([0,\infty)\times \R)$.

\section*{Acknowledgements}
This work was supported by the Royal Society under a Newton International Fellowship (reference number 201625) and by the Conselho Nacional de Desenvolvimento Científico e Tecnológico (process number 308884/2022-1). Finally, the author would like to thank Professor Igor Leite Freire for all the fruitful discussions on the topic.


\begin{thebibliography}{99}


\bibitem{BHP} R. F. Barostichi, A. A. Himonas and G. Petronilho, Autonomous Ovsyannikov theorem and applications to nonlocal evolution equations and systems, \textit{J. Funct. Anal.}, \textbf{v. 270}, 330--358, (2016).

\bibitem{BHP1} R. F. Barostichi, A. A. Himonas and G. Petronilho, Global analyticity for a generalized Camassa–Holm equation and decay of the radius of spatial analyticity, \textit{J. Diff. Eq.}, \textbf{v. 263}, 732--764, (2017).

\bibitem{CH} R. Camassa and D. D. Holm, An integrable shallow water equation with peaked solitons, \textit{Phys. Rev. Lett.}, \textbf{v. 71}, 1661--1664, (1993).

\bibitem{CE} A. Constantin, J. Escher, Global existence and blow-up for a shallow water equation, \textit{Annali Sc. Norm. Sup. Pisa}, \textbf{v. 26}, 303--328, (1998).

\bibitem{Eu} P. L. da Silva, Local well-posedness and global analyticity for solutions of a generalized 0-equation, \textit{Proc. Royal Soc. Edinb.: Sec. A Math.}, (2023). DOI: 10.1017/prm.2022.64

\bibitem{dSF} P. L. da Silva and I. L. Freire, Existence, persistence, and continuation of solutions for a generalized 0-Holm-Staley equation, \textit{J. Diff. Eq.}, \textbf{v. 320}, 371--398, (2022).

\bibitem{DP} A. Degasperis and M. Procesi, Asymptotic Integrability, in Symmetry and Perturbation Theory, edited by A. Degasperis and G.  Gaeta, 23–37, Singapore: World Scientific, (1999).

\bibitem{DHH} A. Degasperis, D. D. Holm and A. N.W. Hone, A new integral equation with peakon solutions, \textit{Theor. Math. Phys.}, \textbf{v. 133}, 1463--1474, (2002).

\bibitem{EY} J. Escher and Z. Yin, Well-posedness, blow-up phenomena, and global solutions for the b-equation, \textit{J. Reine Angw. Math.}, \textbf{v. 624}, 51--80, (2008).

\bibitem{F} I. L. Freire, Remarks on strong global solutions of the $b$-equation, \textit{Appl. Math. Lett.}, \textbf{v. 146}, 108820, (2023).

\bibitem{GLT} G. Gui, Y. Liu, L. Tian, Global existence and blow-up phenomena for the peakon b-family of equations, \textit{Indiana Univ. Math. J.}, \textbf{v. 57}, 1209--1234, (2008).

\bibitem{HM} A. Himonas and G. Misiolek, Analyticity of the Cauchy problem for an integrable evolution equation, \textit{Math. Ann.}, \textbf{v. 327}, 575--584, (2003).

\bibitem{HS1} D. D. Holm and M. F. Staley, Wave structure and nonlinear balances in a family of evolutionary PDEs, \textit{SIAM J. Appl. Dynam. Syst.}, \textbf{v. 2}, 323-380, (2003).

\bibitem{HS2} D. Holm and M. Staley, Nonlinear balance and exchange of stability in dynamics of solitons, peakons, ramp/cliffs and leftons in 1+1 nonlinear evolutionary PDE, \textit{Phys. Lett. A}, \textbf{v. 308}, 437--444, (2003).

\bibitem{KM} T. Kato and K. Masuda, Nonlinear Evolution Equations and Analyticity. I, \textit{Ann. Inst. Henri Poincarè}, \textbf{v. 3}, 455--467, (1986).

\bibitem{LY} Y. Liu and Z. Yin, Global existence and blow-up phenomena for the Degasperis-Procesi equation, \textit{Commun. Math. Phys.}, \textbf{v. 267}, 801--820, (2006).



\end{thebibliography}
\end{document}